
\magnification=1200

\catcode`\À=\active \defÀ{\`A}    \catcode`\à=\active \defà{\`a} 
\catcode`\Â=\active \defÂ{\^A}    \catcode`\â=\active \defâ{\^a} 
\catcode`\Æ=\active \defÆ{\AE}    \catcode`\æ=\active \defæ{\ae}
\catcode`\Ç=\active \defÇ{\c C}   \catcode`\ç=\active \defç{\c c}
\catcode`\È=\active \defÈ{\`E}    \catcode`\è=\active \defè{\`e} 
\catcode`\É=\active \defÉ{\'E}    \catcode`\é=\active \defé{\'e} 
\catcode`\Ê=\active \defÊ{\^E}    \catcode`\ê=\active \defê{\^e} 
\catcode`\Ë=\active \defË{\"E}    \catcode`\ë=\active \defë{\"e} 
\catcode`\Î=\active \defÎ{\^I}    \catcode`\î=\active \defî{\^\i}
\catcode`\Ï=\active \defÏ{\"I}    \catcode`\ï=\active \defï{\"\i}
\catcode`\Ô=\active \defÔ{\^O}    \catcode`\ô=\active \defô{\^o} 
\catcode`\Ù=\active \defÙ{\`U}    \catcode`\ù=\active \defù{\`u} 
\catcode`\Û=\active \defÛ{\^U}    \catcode`\û=\active \defû{\^u} 
\catcode`\Ü=\active \defÜ{\"U}    \catcode`\ü=\active \defü{\"u} 

\catcode`\ =\active \def { }

\hsize=11.25cm    
\vsize=18cm       
\parindent=12pt   \parskip=5pt     

\hoffset=.5cm   
\voffset=.8cm   

\pretolerance=500 \tolerance=1000  \brokenpenalty=5000

\catcode`\@=11

\font\eightrm=cmr8         \font\eighti=cmmi8
\font\eightsy=cmsy8        \font\eightbf=cmbx8
\font\eighttt=cmtt8        \font\eightit=cmti8
\font\eightsl=cmsl8        \font\sixrm=cmr6
\font\sixi=cmmi6           \font\sixsy=cmsy6
\font\sixbf=cmbx6

\font\tengoth=eufm10 
\font\eightgoth=eufm8  
\font\sevengoth=eufm7      
\font\sixgoth=eufm6        \font\fivegoth=eufm5

\skewchar\eighti='177 \skewchar\sixi='177
\skewchar\eightsy='60 \skewchar\sixsy='60

\newfam\gothfam           \newfam\bboardfam

\def\tenpoint{
  \textfont0=\tenrm \scriptfont0=\sevenrm \scriptscriptfont0=\fiverm
  \def\rm{\fam\z@\tenrm}
  \textfont1=\teni  \scriptfont1=\seveni  \scriptscriptfont1=\fivei
  \def\oldstyle{\fam\@ne\teni}\let\old=\oldstyle
  \textfont2=\tensy \scriptfont2=\sevensy \scriptscriptfont2=\fivesy
  \textfont\gothfam=\tengoth \scriptfont\gothfam=\sevengoth
  \scriptscriptfont\gothfam=\fivegoth
  \def\goth{\fam\gothfam\tengoth}
  
  \textfont\itfam=\tenit
  \def\it{\fam\itfam\tenit}
  \textfont\slfam=\tensl
  \def\sl{\fam\slfam\tensl}
  \textfont\bffam=\tenbf \scriptfont\bffam=\sevenbf
  \scriptscriptfont\bffam=\fivebf
  \def\bf{\fam\bffam\tenbf}
  \textfont\ttfam=\tentt
  \def\tt{\fam\ttfam\tentt}
  \abovedisplayskip=12pt plus 3pt minus 9pt
  \belowdisplayskip=\abovedisplayskip
  \abovedisplayshortskip=0pt plus 3pt
  \belowdisplayshortskip=4pt plus 3pt 
  \smallskipamount=3pt plus 1pt minus 1pt
  \medskipamount=6pt plus 2pt minus 2pt
  \bigskipamount=12pt plus 4pt minus 4pt
  \normalbaselineskip=12pt
  \setbox\strutbox=\hbox{\vrule height8.5pt depth3.5pt width0pt}
  \let\bigf@nt=\tenrm       \let\smallf@nt=\sevenrm
  \normalbaselines\rm}

\def\eightpoint{
  \textfont0=\eightrm \scriptfont0=\sixrm \scriptscriptfont0=\fiverm
  \def\rm{\fam\z@\eightrm}
  \textfont1=\eighti  \scriptfont1=\sixi  \scriptscriptfont1=\fivei
  \def\oldstyle{\fam\@ne\eighti}\let\old=\oldstyle
  \textfont2=\eightsy \scriptfont2=\sixsy \scriptscriptfont2=\fivesy
  \textfont\gothfam=\eightgoth \scriptfont\gothfam=\sixgoth
  \scriptscriptfont\gothfam=\fivegoth
  \def\goth{\fam\gothfam\eightgoth}
  
  \textfont\itfam=\eightit
  \def\it{\fam\itfam\eightit}
  \textfont\slfam=\eightsl
  \def\sl{\fam\slfam\eightsl}
  \textfont\bffam=\eightbf \scriptfont\bffam=\sixbf
  \scriptscriptfont\bffam=\fivebf
  \def\bf{\fam\bffam\eightbf}
  \textfont\ttfam=\eighttt
  \def\tt{\fam\ttfam\eighttt}
  \abovedisplayskip=9pt plus 3pt minus 9pt
  \belowdisplayskip=\abovedisplayskip
  \abovedisplayshortskip=0pt plus 3pt
  \belowdisplayshortskip=3pt plus 3pt 
  \smallskipamount=2pt plus 1pt minus 1pt
  \medskipamount=4pt plus 2pt minus 1pt
  \bigskipamount=9pt plus 3pt minus 3pt
  \normalbaselineskip=9pt
  \setbox\strutbox=\hbox{\vrule height7pt depth2pt width0pt}
  \let\bigf@nt=\eightrm     \let\smallf@nt=\sixrm
  \normalbaselines\rm}

\tenpoint

\def\pc#1{\bigf@nt#1\smallf@nt}         \def\pd#1 {{\pc#1} }

\catcode`\;=\active
\def;{\relax\ifhmode\ifdim\lastskip>\z@\unskip\fi
\kern\fontdimen2  -1.2 \fontdimen3 \string;}

\catcode`\:=\active
\def:{\relax\ifhmode\ifdim\lastskip>\z@\unskip\fi\penalty\@M\ \fi\string:}

\catcode`\!=\active
\def!{\relax\ifhmode\ifdim\lastskip>\z@
\unskip\fi\kern\fontdimen2  -1.1 \fontdimen3 \string!}

\catcode`\?=\active
\def?{\relax\ifhmode\ifdim\lastskip>\z@
\unskip\fi\kern\fontdimen2  -1.1 \fontdimen3 \string?}

\catcode`\«=\active 
\def«{\raise.4ex\hbox{%
 $\scriptscriptstyle\langle\!\langle$}}

\catcode`\»=\active 
\def»{\raise.4ex\hbox{%
 $\scriptscriptstyle\rangle\!\rangle$}}

\frenchspacing

\def\raggedbottom{\topskip 10pt plus 36pt\r@ggedbottomtrue}

\def\pointir{\unskip . --- \ignorespaces}

\def\Medbreak{\vskip-\lastskip\medbreak}

\long\def\th#1 #2\enonce#3\endth{
   \Medbreak\noindent
   {\pc#1} {#2\unskip}\pointir{\it #3}\smallskip}

\def\proof{\vskip-\lastskip\smallskip\noindent
 {\it Proof} : }

\def\decale#1{\smallbreak\hskip 28pt\llap{#1}\kern 5pt}
\def\decaledecale#1{\smallbreak\hskip 34pt\llap{#1}\kern 5pt}
\def\puce{\smallbreak\hskip 6pt{$\scriptstyle\bullet$}\kern 5pt}

\def\eqalign#1{\null\,\vcenter{\openup\jot\m@th\ialign{
\strut\hfil$\displaystyle{##}$&$\displaystyle{{}##}$\hfil
&&\quad\strut\hfil$\displaystyle{##}$&$\displaystyle{{}##}$\hfil
\crcr#1\crcr}}\,}

\newcount\numeroderemarque
\def\remarque{\advance\numeroderemarque by1\smallbreak
{\it Remarque\/}\ \number\numeroderemarque~:}

\newcount\numeroderemarque
\def\remark{\advance\numeroderemarque by1\smallbreak
{\it Remark\/}\ \number\numeroderemarque~:}

\catcode`\@=12

\showboxbreadth=-1  \showboxdepth=-1

\newcount\numerodesection \numerodesection=1
\def\section#1{\bigbreak
 {\bf\number\numerodesection.\ \ #1}\nobreak\medskip
 \advance\numerodesection by1}

\mathcode`A="7041 \mathcode`B="7042 \mathcode`C="7043 \mathcode`D="7044
\mathcode`E="7045 \mathcode`F="7046 \mathcode`G="7047 \mathcode`H="7048
\mathcode`I="7049 \mathcode`J="704A \mathcode`K="704B \mathcode`L="704C
\mathcode`M="704D \mathcode`N="704E \mathcode`O="704F \mathcode`P="7050
\mathcode`Q="7051 \mathcode`R="7052 \mathcode`S="7053 \mathcode`T="7054
\mathcode`U="7055 \mathcode`V="7056 \mathcode`W="7057 \mathcode`X="7058
\mathcode`Y="7059 \mathcode`Z="705A


\catcode`\À=\active \defÀ{\`A}    \catcode`\à=\active \defà{\`a} 
\catcode`\Â=\active \defÂ{\^A}    \catcode`\â=\active \defâ{\^a} 
\catcode`\Æ=\active \defÆ{\AE}    \catcode`\æ=\active \defæ{\ae}
\catcode`\Ç=\active \defÇ{\c C}   \catcode`\ç=\active \defç{\c c}
\catcode`\È=\active \defÈ{\`E}    \catcode`\è=\active \defè{\`e} 
\catcode`\É=\active \defÉ{\'E}    \catcode`\é=\active \defé{\'e} 
\catcode`\Ê=\active \defÊ{\^E}    \catcode`\ê=\active \defê{\^e} 
\catcode`\Ë=\active \defË{\"E}    \catcode`\ë=\active \defë{\"e} 
\catcode`\Î=\active \defÎ{\^I}    \catcode`\î=\active \defî{\^\i}
\catcode`\Ï=\active \defÏ{\"I}    \catcode`\ï=\active \defï{\"\i}
\catcode`\Ô=\active \defÔ{\^O}    \catcode`\ô=\active \defô{\^o} 
\catcode`\Ù=\active \defÙ{\`U}    \catcode`\ù=\active \defù{\`u} 
\catcode`\Û=\active \defÛ{\^U}    \catcode`\û=\active \defû{\^u} 
\catcode`\Ü=\active \defÜ{\"U}    \catcode`\ü=\active \defü{\"u} 

\def\ogoth{{\goth o}}
\def\units{{\ogoth^\times}}

\def\pgoth{{\goth p}}

\def\Q{{\bf Q}}

\def\Z{{\bf Z}}

\def\F{{\bf F}}

\def\Gal{\mathop{\rm Gal}\nolimits}
\def\Ker{\mathop{\rm Ker}\nolimits}

\def\ogoth{{\goth o}}
\def\units{{\ogoth^\times}}

\def\pgoth{{\goth p}}

\def\Q{{\bf Q}}

\def\P{{\bf P}}

\def\Z{{\bf Z}}

\def\F{{\bf F}}

\def\Gal{\mathop{\rm Gal}\nolimits}
\def\Pic{\mathop{\rm Pic}\nolimits}
\def\Ker{\mathop{\rm Ker}\nolimits}

\def\zero{\{0\}}

\def\ei{e_1}
\def\eii{e_2}
\def\zmodii{\Z / 2\Z}
\def\zmodiixzmodii{(\zmodii)^2}
\def\mod{\mathop{\rm mod.}\nolimits}
\def\pmod#1{\;(\mod#1)}
\def\mmod#1{\;(\mod^{\!\!\times}{\pgoth^{#1}})}

\newcount\refno 
\long\def\ref#1:#2<#3>{                                        
\global\advance\refno by1\par\noindent                              
\llap{[{\bf\number\refno}]\ }{#1} \pointir{\it #2} #3\goodbreak }

\def\citer#1(#2){[{\bf\number#1}\if#2\empty\relax\else,\ #2\fi]}

\newcount\formuleno
\def\numeroter{\global\advance\formuleno by1
 \leqno{(\oldstyle\number\formuleno)}}
\def\formule#1{$(\oldstyle\number#1)$}

\newbox\bibbox
\setbox\bibbox\vbox{\bigbreak
\centerline{{\pc BIBLIOGRAPHIC} {\pc REFERENCES}}

\ref{\pc COLLIOT}-{\pc TH{\'E}L{\`E}NE} (J-L.) and {\pc CORAY} (D. F.):
L'{\'e}quivalence rationnelle sur les points ferm{\'e}s des surfaces
rationnelles fibr{\'e}es en coniques,
<Compositio Math., {\bf 39} (3), 1979, pp.~301--332>
\newcount\ctcoray  \global\ctcoray=\refno

\ref {\pc COLLIOT}-{\pc THÉLÈNE} (J-L.) and {\pc ISCHEBECK} (F.):
L'équivalence rationnelle sur les cycles de
dimension zéro des variétés algébriques réelles,
<Comptes-rendus de l'Acad.\ des sci.\ {\bf 292} (1981), 723--725.>
\newcount\ctischebeck \global\ctischebeck=\refno

\ref{\pc COLLIOT}-{\pc TH{\'E}L{\`E}NE} (J-L.) and {\pc SANSUC} (J-J.):
La descente sur les vari{\'e}t{\'e}s rationnelles,
<dans Journ{\'e}es de G{\'e}om{\'e}trie alg{\'e}brique d'Anger, Alpen
aan den Rijn~: Sijthoff \& Noordhoff, 1980.>
\newcount\ctsansucdescente  \global\ctsansucdescente=\refno

\ref{\pc COLLIOT}-{\pc TH{\'E}L{\`E}NE} (J-L.) and {\pc SANSUC} (J-J.):
On the {C}how groups of certain rational surfaces: a sequel to
a paper of {S}.\ {B}loch,
<Duke Math.\ J., {\bf 48} (2), 1981, pp.~421--447>
\newcount\ctsansucsequel  \global\ctsansucsequel=\refno

\ref{\pc COLLIOT}-{\pc TH{\'E}L{\`E}NE} (J-L.), {\pc SANSUC} (J-J.)
and {\pc SWINNERTON-\pc DYER} (H. P. F.):
Intersections of two quadrics and Châtelet surfaces. I, II.
<J.\ Reine Angew.\ Math.\ 373 (1987), pp.~37--107; {\it ibid.} 374 (1987),
pp.~72--168.>
\newcount\ctssd \global\ctssd=\refno

\ref{\pc CORAY} (D. F.) and {\pc TSFASMAN} (M. A.):
Arithmetic on singular del {P}ezzo surfaces,
<Proc.\ London Math.\ Soc.\ (3) {\bf 57} (1), 1988, pp.~25--87.>
\newcount\coraytsfasman  \global\coraytsfasman=\refno

\ref{\pc DALAWAT} (C. S.):
Le groupe de Chow d'une surface de Ch\^atelet sur un corps local,
<Indagationes math., {\bf 11} (2), 2000, pp.~171--185.>
\newcount\chatelet  \global\chatelet=\refno

\ref{\pc SWINNERTON-\pc DYER} (H. P. F.):
Two descent from Fermat to now, <Mathematisches Institut
Universit{\"a}t G{\"o}ttingen, Summer Term 2004, pp.~95--102,
Universit{\"a}tsdrucke G{\"o}ttingen, G{\"o}ttingen, 2004.>
\newcount\swinn  \global\swinn=\refno

\ref{\pc SANSUC} (J-J.):
{\`A} propos d'une conjecture arithm{\'e}tique sur le groupe de
{C}how d'une surface rationnelle,
<S{\'e}minaire de Th{\'e}orie des nombres de Bordeaux, Expos{\'e} 33, 1982.>
\newcount\sansuc  \global\sansuc=\refno

} 

\centerline{\bf The Chow group of a Ch\^atelet surface}
\smallskip
\centerline{\bf over a number field}
\bigskip
\centerline{Chandan Singh Dalawat}
\bigskip
{\eightpoint We compute the Chow group of a Ch{\^a}telet surface over
a dyadic field.  Combined with the previous work of Bloch,
Colliot-Th{\'e}l{\`e}ne, Coray, Ischebeck, Sansuc, Swinnerton-Dyer,
and the author, this allows one to compute the Chow group of any
Ch{\^a}telet surface over any number field.}

\bigskip
{\bf 1.  Introduction}

$2$-descent (also called a first descent) on elliptic curves $E$
defined over $\Q$ which have the form $y^2=(x-c_1)(x-c_2)(x-c_3)$
(where the $c_i$ are distinct), is a classical theme which goes back
to Pierre Fermat, with a major contribution by John Tate.  If carried
out successfully, one ends up computing the finitely generated
commutative group $E(\Q)$ of rational points on $E$.

There is a simpler version of $2$-descent, applicable to Ch{\^a}telet
surfaces, i.e.~smooth proper surfaces $X$ birational to the affine
surface  defined by 
$$
y^2-dz^2=(x-c_1)(x-c_2)(x-c_3)\qquad(d\in\Q^\times).
\numeroter \newcount\chatsurf \global\chatsurf=\formuleno
$$
This version has also been pursued by various authors, whose
contributions will be recalled in the course of this Note.  We intend
to carry it out successfully, computing thereby the Chow group of $X$,
the finite $\F_2$-space $A_0(X)_0$ of degree-$0$ $0$-cycles modulo
rational equivalence.  

Of the real $2$-descent, the one for elliptic curves, Peter
Swinnerton-Dyer wrote recently that ``the statements of the theory
over an arbitrary number field are not very different, except that the
analogues of certain explicit results relating to the prime~$2$ are
not known'' \citer\swinn().  

The main contribution of this Note is to obtain those explicit results
at the prime~$2$ for $2$-descent on Ch{\^a}telet surfaces~; they allow
us to compute the Chow group of any Ch{\^a}telet surface over any
number field.

\medskip
{\bf 2. Statement of the local results}
\medskip
Let $K$ be a finite extension of the field $\Q_p$ ($p$ prime) of
$p$-adic numbers, $\ogoth$ the ring of integers of $K$ (i.e.~the
integral closure of $\Z_p$ in $K$), $\pgoth$ the unique maximal ideal
of $\ogoth$, and $k=\ogoth/\pgoth$ the residue field of $K$.  Denote
by $v:K^\times\to\Z$ the surjective valuation of $K$.

Given $d\in K^\times$ and three distinct numbers $c_1,c_2,c_3\in K$,
we get a smooth affine $K$-surface
$$
y^2-dz^2=(x-c_1)(x-c_2)(x-c_3)
\numeroter \newcount\locchatsurf \global\locchatsurf=\formuleno
$$
which is $K$-birational to $\P_2$ if $d\in K^{\times2}$~; if
$d\notin K^{\times2}$, it is $K(\sqrt{d})$-birational to $\P_2$, and
the group $A_0(X)_0$ is killed by~$2$.

Let $X$ be any smooth projective $K$-surface which is $K$-birational
to \formule\locchatsurf.  Thus the Chow group $A_0(X)_0$ of
$0$-cycles of degree~$0$ modulo rational equivalence, which we are
interested in computing, is trivial if $d\in K^{\times2}$~; we
therefore assume that $d\notin K^{\times2}$.

To fix ideas, we take $X$ to be the surface
defined in $\P({\cal O}(2)\oplus{\cal O}(2)\oplus{\cal O})$
(coordinates $y:z:t$) over the projective $K$-line $\P_{1,K}$
(coordinates $x:x'$) by the equation
$$
y^2-dz^2=(x-c_1x')(x-c_2x')(x-c_3x')t^2;
\numeroter \newcount\locchatsurfii \global\locchatsurfii=\formuleno
$$
this model was constructed by Colliot-Th{\'e}l{\`e}ne and Sansuc in
\citer\ctsansucdescente().  This choice is immaterial for our purposes, as
we are only interested in computing the group $A_0(X)_0$ which is
independent of the choice of the projective model, by a result of
Colliot-Th{\'e}l{\`e}ne and Coray \citer\ctcoray(). The surface
\formule\locchatsurfii\ comes equipped with a morphism to $\P_1$ whose
fibres are conics~; there are four degenerate fibres, namely the ones
above $c_1$, $c_2$, $c_3$ and $\infty$.

The change of variables $x_1=x-c_1$ allows us to write
\formule\locchatsurf\ as 
$$
y^2-dz^2=x_1(x_1-c_1^\prime)(x_1-c_2^\prime).
$$ 
Upto permuting $c_1^\prime$, $c_2^\prime$, we may suppose that
$v(c_2^\prime)\ge v(c_1^\prime)$.  Moreover, if
$v(c_2^\prime)>v(c_1^\prime)$, the change of variables
$x_2=x_1-c_1^\prime$ transforms this equation into
$$
y^2-dz^2=x_2(x_2-c_1'')(x_2-c_2'')
$$ 
in which $v(c_2'')=v(c_1'')$.  In other words, we can suppose without
any loss of generality that we have $c_1=0$ and $v(c_2)=v(c_3)$ in
\formule\locchatsurf.  Denoting this common valuation by $r$, we shall
henceforth work with
$$
y^2-dz^2=x(x-e_1)(x-e_2)\qquad (r=v(e_i)).
\numeroter \newcount\locchatsurf \global\locchatsurf=\formuleno
$$

With these conventions, let us recall the cases in which the group
$A_0(X)_0$ has been computed.

\th PROPOSITION 1 (\citer\coraytsfasman(Prop.~4.7) for $p$~odd,
\citer\chatelet(Prop.~1) for $p=2$)
\enonce
Suppose that the extension $K(\sqrt{d})$ is unramified. The group
$A_0(X)_0$ is then isomorphic 
$$\vbox{\halign{\hfil\it #\/\rm)\ &{to}#&\qquad\hss$#$\hss\qquad
&{\it if\/ #}\hfil\cr
  i&&0^{\phantom{2}} & $r$ is even and\/ $v(\ei-\eii)=r$,\cr
 ii&&\zmodii^{\phantom{2}} & $r$ is even and\/ $v(\ei-\eii)>r$,\cr
iii&&\zmodiixzmodii & $r$ is odd.\cr}}
$$
\endth

\th PROPOSITION 2 (\citer\chatelet(Prop.~2))
\enonce
Suppose that $p$ is odd and that the extension $L=K(\sqrt{d})$ is
ramified. The group $A_0(X)_0$ is then isomorphic
$$
\vbox{\halign{\hfil\it #\/\rm)\ &{to}#&\qquad\hss$#$\hss\qquad&
{\it if\/ #}\hfil\cr
  i&&\zmodii^{\phantom{2}} & $\ei/\eii\equiv1\pmod\pgoth$ 
                                and\/ $\ei\in N_{L|K}(L^\times)$,\cr 
 ii&&\zmodiixzmodii & $\ei/\eii\equiv1\pmod\pgoth$ 
                                and\/ $\ei\notin N_{L|K}(L^\times)$,\cr 
iii&&\zmodiixzmodii & $\ei/\eii\not\equiv1\pmod\pgoth$.\cr}}
$$
\endth

Suppose that $p=2$ and that the extension $L=K(\sqrt{d})$ is ramified.
Let $\chi:K^\times\rightarrow\Z/2\Z$ be the homomorphism whose kernel
is the group of norms from $L^\times$.  The restriction of $\chi$ to
the units $\ogoth^\times$ is $\neq0$, because $L|K$ is ramified.
Since $k^\times$ is of odd order, $\chi|_{\ogoth^\times}$ factors via
the group $U_1$ of $1$-units, and indeed via $U_1/U_{n+1}$, but not via
$U_1/U_n$, for a suitable $n>0$.  Here $U_n$ denotes the group of
units which are $\equiv1\pmod{\pgoth^n}$.

As the archimedean local fields were treated in \citer\ctischebeck(),
the missing ingredient is provided by the next proposition, which is
our main result.  This proposition not only completes the
determination of the Chow group of a Ch{\^a}tlet surface over a local
field, but, becuase of a local-to-global principle which we will
recall below, also over number fields.

\th PROPOSITION 3
\enonce
Suppose that\/ $p=2$ and that the extension\/ $L=K(\sqrt{d})$ is
ramified. Suppose that\/ $\chi$ factors via\/ $U_1/U_{n+1}$, but not via\/
$U_1/U_n$.  The group\/ $A_0(X)_0$ is then isomorphic
$$
\vbox{\halign{\hfil\it #\/\rm)\ &{to}#&\qquad\hss$#$\hss\qquad&
{\it if\/ #}\hfil\cr
  i&&\zmodii^{\phantom{2}} & $\ei/\eii\equiv1\pmod{\pgoth^{2n+1}}$ 
                                and\/ $\ei\in N_{L|K}(L^\times)$,\cr 
 ii&&\zmodiixzmodii & $\ei/\eii\equiv1\pmod{\pgoth^{2n+1}}$ 
                                and\/ $\ei\notin N_{L|K}(L^\times)$,\cr 
iii&&\zmodiixzmodii & $\ei/\eii\not\equiv1\pmod{\pgoth^{2n+1}}$.\cr}}
$$
\endth
As we shall see in the course of the proof, these three propositions,
as well as the result at the places at $\infty$, actually compute
$A_0(X)_0$ as a subgroup of $H^1(K,S(\bar K))$, where $\bar K$ is an
algebraic closure of $K$ and $S$ is the $K$-torus whose group of
characters is the $\Gal(\bar K|K)$-module $\Pic(\bar X)$, with $\bar
X=X\times_K\bar K$.  This is crucial for the application to computing
the Chow group of a Ch{\^a}telet surface over a number field.

\medbreak
{\bf 3. The method of computation}
\medskip

It relies on the work of Colliot-Th{\'e}l{\`e}ne and Coray \citer\ctcoray()
and Colliot-Th{\'e}l{\`e}ne and Sansuc \citer\ctsansucdescente(),
\citer\ctsansucsequel(); it has been explained at length in
\citer\chatelet().  As this last paper was written in a language which is
poorly understood in many parts of the world, ce qui n'est pas sans
rapport avec les p{\'e}rip{\'e}ties de l'auteur, we shall gave a brief
summary here.

One can replace $X$ by any smooth proper $K$-surface $K$-birational to
$X$~\citer\ctcoray(Prop.~6.3), which justifies the choice of the
particular model $X$ (cf.~\formule\locchatsurfii) of the equation
\formule\locchatsurf\ that we have made.  

Next, denoting by $O$ the singular point of the fibre at infinity of
the conic bundle $f:X\to\P_1$, the map
$$
\gamma:X(K)\to A_0(X)_0,\quad \gamma(Q)=Q-O
$$
is surjective \citer\ctcoray(Th{\'e}or{\`e}me~C).  ``The
characteristic homomorphism'' is a natural injection 
$$
\varphi:A_0(X)_0\to H^1(K,S(\bar K))
$$ 
\citer\ctsansucsequel(n$^{\rm o}$~IV).  With the identifications 
$$
\iota:H^1(K,S(\bar K))\to(K^\times\!/N_{L|K}(L^\times))^2\to(\Z/2\Z)^2,
$$
the composite map $X(K)\to(\Z/2\Z)^2$ is given by
$$
\def\classe#1{\chi(#1)}
(y:z:t;x)\mapsto\cases{
(\classe{1},\, \classe{1}) & if $x=\infty$,\cr
(\classe{\ei\eii},\, \classe{-\ei}) & if $x=0$,\cr
(\classe{\ei},\, \classe{\ei(\ei-\eii})& if $x=\ei$,\cr
(\classe{x},\, \classe{x-\ei}) & otherwise.\cr}\leqno{(1)}
$$
\citer\ctsansucdescente(n$^{\rm o}$~IV).  As all the points in the
same fibre of the map $f:X(K)\to\P_1(K)$ are mutually equivalent
$0$-cycles, what we have to compute is the image of the induced map
$[\phantom{x}]:f(X(K))\to(\Z/2\Z)^2$.  The subset
$f(X(K))\subset\P_1(K)$ consists of $\infty$, $0$, $\ei$, $\eii$ and
all those $x\in K$, different from $0,$ $\ei$, $\eii$, for which
$\chi(x(x-\ei)(x-\eii))=0$.  As explained in \citer\chatelet(), what
we have to compute is the subgroup generated by the image $[f(X(K))]$.

All this is valid for $K$ any extension of $\Q_p$, and $L=K(\sqrt{d})$
any quadratic extension of $K$.  As explained in
\citer\chatelet(Remarque~5), when $L$ is a ramified extension, we can
suppose that $v(e_i)=0$.  Further, when $p=2$, we can assume that
$e_i\in U_1$, as every element of $k^\times$ is in
$N_{L|K}(L^\times)$. 

We assume from now onwards that $p=2$, that $L|K$ is ramified, and
that the $e_i$ are (distinct) $1$-units of $K$.

\th LEMMA 1
\enonce
Suppose that $\chi$ factors via
$(\ogoth/\pgoth^n)^\times\rightarrow\Z/2\Z$.  Then, for every $x\in
f(X(K))$, we have
$$
[x] =
\cases{
(0,0) & if $v(x)\le-n,$\cr 
\phantom{0}[0] & if $v(x)\ge n$.\cr 
}$$
\endth
\proof\ For $x\in K$ distinct from $0,\ei,\eii$, we have 
$$
\chi(x(x-\ei)(x-\eii))=\chi(x)+\chi(x-\ei)+\chi(x-\eii);
$$
if further $x\in f(X(K))$, this sum is $0$.  If moreover $v(x)\le-n$, then
each term is $0$, since $\chi(x-e_i)=\chi(x)$ in this case, so
$[x]=(0,0)$.  If however $v(x)\ge n$, we have $\chi(x-e_i)=\chi(-e_i)$
and therefore $\chi(x)=\chi(\ei\eii)$; accordingly, $[x]=[0]$.

\th PROPOSITION 4
\enonce
Suppose that $\chi$ factors via\/ $U_1/U_{n+1}$ but does not factor
via\/ $U_1/U_n$.  Then the subgroup generated by\/ $[f(X(K))]$ is\/~:
$$
\vbox{\halign{\hfil\it #\/\rm)\ &\qquad\hss$#$\hss\qquad&
{\it if\/ #}\hfil\cr
  i&\zero\times\zmodii & $\ei\equiv\eii\pmod{\pgoth^{2n+1}}$ and\/
  $\chi(\ei)=0$,\cr  
 ii&\zmodiixzmodii & $\ei\equiv\eii\pmod{\pgoth^{2n+1}}$ and\/
  $\chi(\ei)\neq0$,\cr  
iii&\zmodiixzmodii & $\ei\not\equiv\eii\pmod{\pgoth^{2n+1}}$.\cr}}
$$
\endth\noindent
We shall indentify $k^\times$ with the prime-to-$2$ torsion subgroup
of $\units$.

Putting $H=\Ker(\chi)$, notice that $1+t\bar\pi^n\notin H$ for some
$t\in k^\times$, since otherwise $\chi$ would factor via $U_1/U_n$.
Fix such a $t\in k^\times$~; notice that $1-t\bar\pi^n\notin H$, for
$1-t\bar\pi^n\equiv1+t\bar\pi^n\;(\mod\pgoth^{n+1})$, as
$2\equiv0\;(\mod\pgoth^e)$.  If $a\in H$, then $a(1+t\bar\pi^n)\notin
H$, for their quotient is $1+t\bar\pi^n$.  Finally, $H\rightarrow
U_1/U_n$ is a surjection, with kernel $H\cap(U_n/U_{n+1})$ of order
$2^{f-1}$, where $q=2^f$ is the number of elements in the residue
field $k$.

Put $x=t^{-1}\pi^{-n}$.  Observe that $x\in f(X(K))$ and $[x]=(0,1)$ in all
cases~: we have $\chi(x)=0$, and $\chi(x-e_i)=\chi(1-te_i\pi^n)=1$, because
$1-te_i\pi^n\equiv1-t\pi^n\;(\mod\pgoth^{n+1})$.

{\it Proof of i)\/} We have $[0]=(0,*)$, $[\ei]=(0,*)$,
$[\eii]=(0,*)$; let us show that $(1,*)$ is not (i.e.~$(1,0)$ and
$(1,1)$ are not) in the image $[f(X(K))]$.  It is sufficient to show that if
$x\in K$ is distinct from $0$, $\ei$, $\eii$, if its valuation is
between $-n$ and $n$ (cf.~lemma~1), and if $\chi(x)=1$, then $x\notin
f(X(K))$.  Recalling that $v(\ei-\eii)>2n$ and that $v(e_i)=0$, we deduce
that $v(x-\ei)=v(x-\eii)$ is given by
$$
\vbox{\halign{&\hfil$#$\hfil\quad\cr
v(x)&-n&\cdots&-1&0&1&\cdots&n\cr
\noalign{\vskip-5pt}
\multispan8\hrulefill.\cr
v(x-e_i)&-n&\cdots&-1&\ 0, 1,\cdots,n\ &0&\cdots&0\cr
}}
$$
Only the case $v(x)=0$ needs some explanation~: the fact that
$\chi(x)=1$ whereas $\chi(e_i)=0$ means that $x\not\equiv
e_i\mod\pgoth^{n+1}$, i.e.~$v(x-e_i)<n+1$.

Writing $r=v(x-e_i)$, we have
$$
\pi^{-r}(x-\eii)=\pi^{-r}(x-\ei)+\pi^{-r}(\ei-\eii)
$$
with $v(\pi^{-r}(\ei-\eii))>n$.  So we get $x-\ei\equiv
x-\eii\mmod{n+1}$, which implies that $\chi(x-\ei)=\chi(x-\eii)$, and
hence $x\notin f(X(K))$.  This failure contains the seeds of our success in
case {\it iii\/}) below.

{\it Proof of ii)\/}  We have $[\ei]=(1,*)$, as $\chi(\ei)=1$. 

{\it Proof of iii)\/} It remains to show that $(1,*)$ is (i.e.~$(1,0)$
or $(1,1)$ is) in $[f(X(K))]$, which is clearly the case if $\chi(\ei)=1$,
and also if $\chi(\eii)=1$.  Suppose then that $\chi(e_i)=0$.

Write $\ei-\eii=s_{1,2}u_{1,2}\pi^a$ ($s_{1,2}\in k^\times$,
$u_{1,2}\in U_1$, $1\le a\le2n$) and let us first deal with the case
$a=2n$.  Identify the $\F_2$-space $k$ with the subgroup $U_n/U_{n+1}$
of $U_1/U_{n+1}$ by the map $s\mapsto1+s\pi^n$, and let $h\subset k$
be the codimension-$1$ subspace such that $\chi(1+s\pi^n)=0$ if and
only if $s\in h$, i.e.~$h=H\cap k$.  There is a $t\notin h$ such that
$s_{1,2}t^{-1}\notin h$, for, as $t\notin h$ varies, we get
$2^f-2^{f-1}=2^{f-1}$ distinct elements $s_{1,2}t^{-1}$ of $k^\times$,
of which at most $2^{f-1}-1$ can belong to $h$~: we are saved by the
skin of our teeth.  Fix such a $t$, and put $x=\ei+t\pi^n$, so that
$\chi(x)=1$ --- this is the only place where we are using the fact
that $\chi(\ei)=0$ --- and $\chi(x-\ei)=0$.  Also, $v(x-\eii)=n$ and,
writing $x-\eii=(x-\ei)+(\ei-\eii)$, we get
$$
{x-\eii\over\pi^{n}}\equiv t(1+s_{1,2}t^{-1}\pi^n)\mod\pgoth^{n+1},
$$
so $\chi(x-\eii)=1$.  Thus, $x\in f(X(K))$ and $[x]=(1,0)$.

Assume now that the valuation $a$ of $\ei-\eii$ is $<2n$ and let $t\in
k^\times$ be such that $1+t\bar\pi^n\notin H$.  Put $s=s_{1,2}t^{-1}$
and take $x=\ei+su\pi^{a-n}$, where $u\in U_1$ is such that
$\chi(x)=1$~; this can always be arranged, at the cost of replacing
$x$ by $x(1+t\pi^n)$, i.e.~$u$ by the $w$ satisfying
$x(1+t\pi^n)=\ei+sw\pi^{a-n}$ ($w\in U_1$; in fact,
$w=u+tu\pi^n+s^{-1}t\pi^{2n-a}$).  We have $v(x-e_i)=a-n$, and
$$
{x-\eii\over\pi^{a-n}}=su+s_{1,2}u_{1,2}\pi^n
\equiv su(1+t\pi^n)\mod\pgoth^{n+1},
$$
so $\chi(x-\eii)+\chi(x-\ei)=1$, i.e.~$x\in f(X(K))$.  We have $[x]=(1,*)$.

This completes the proof of Prop.~4, and thereby also Prop.~3.
\medskip
{\bf 4.  Consequence}
\medskip

The main application of Prop.~3 is to the computation of the Chow
group of $0$-cycles of degree~$0$ on a Ch{\^a}telet surface defined
over a number field, an application made possible by the
local-to-global principle as conjectured in \citer\ctsansucsequel()
and proved in \citer\ctssd().  As this process has been clearly
explained in \citer\sansuc() and
\citer\ctsansucsequel(), we content ourselves with a practical remark.   

Changing notation, let $K$ be a finite extension of $\Q$, $d\in
K^\times$ and let $c_1, c_2, c_3\in K$ be distinct.  Let $X$ be any
smooth proper surface $K$-birational to \formule\locchatsurf.  For
each place $v$ of $K$, we have a map $A_0(X)_0\to A_0(X_v)_0$, where
$X_v=X\times K_v$ and $K_v$ is the completion of $K$ at the place
$v$. We also have a map $j_v:A_0(X_v)_0\to(\Z/2\Z)^3$ (restoring
symmetry) whose image lies in the subgroup $b_1+b_2+b_3=0$.  If $d\in
K_v^{\times2}$, then $A_0(X_v)_0=\zero$ \citer\ctcoray(Prop.~4.7), as
$X_v$ is then $K_v$-birational to $\P_2$.  Assume that $d\notin
K_v^{\times2}$.  

If $v$ is a real place, the results of
\citer\ctischebeck() allow us to compute $j_v$.  Assume that
$v$ does not lie above $\infty$.

Suppose that the extension $K_v(\sqrt{d})$ is unramified.  If $v$ lies
above an odd prime, the results of
\citer\coraytsfasman(Prop.~4.7) (cf.~Prop.~1; see also
\citer\chatelet(n$^{\rm o}$~4)) allow us to compute
$j_v$.  If $v$ lies above $2$, then one applies
\citer\chatelet(n$^{\rm o}$~4).

Suppose now that the extension $K_v(\sqrt{d})$ is ramified.  If $v$ lies
above an odd prime, one can use \citer\chatelet(Prop.~2)
(cf.~Prop.~2).  Finally, if $v$ divides~$2$, one applies Prop.~3 to
compute $j_v$.

By the local-to-global principle as conjectured in \citer\ctsansucsequel()
and proved in \citer\ctssd(), the group $A_0(X)_0$ is the kernel of
the map $\oplus_vj_v$ form $\oplus_vA_0(X_v)_0$ into $(\Z/2\Z)^3$.

\unvbox\bibbox 
\vskip2cm
{\obeylines\parskip=0pt\parindent=0pt
Chandan Singh Dalawat
Harish-Chandra Research Institute
Chhatnag Road, Jhunsi
{\pc ALLAHABAD} 211\thinspace019, India
\vskip5pt
\tt dalawat@mri.ernet.in}

\bye